\documentclass[10pt,twoside,a4paper,reqno]{amsart}
\usepackage{amscd,amsmath,amsthm,amsfonts,latexsym,amssymb,epsfig}
\usepackage{graphics,graphicx}

\theoremstyle{plain}
\newtheorem{theo+}           {Theorem}
\newtheorem{prop+}           {Proposition}
\newtheorem{coro+}           {Corollary}
\newtheorem{lemm+}           {Lemma}
\newtheorem{conjecture}      {Conjecture}

\theoremstyle{definition}
\newtheorem{rema+}           {Remark}
\newtheorem{defi+}           {Definition}
\newtheorem{problem}         {Problem}

\newenvironment{theorem}{\begin{theo+}}{\end{theo+}}
\newenvironment{proposition}{\begin{prop+}}{\end{prop+}}
\newenvironment{corollary}{\begin{coro+}}{\end{coro+}}
\newenvironment{lemma}{\begin{lemm+}}{\end{lemm+}}
\newenvironment{remark}{\begin{rema+}}{\end{rema+}}
\newenvironment{definition}{\begin{defi+}}{\end{defi+}}

\newcommand{\D}{\mathcal{D}}

\newcommand{\calH}{\mathcal{H}}
\newcommand{\RD}{\mathcal{RD}}
\newcommand{\CG}{\mathcal{G}}
\newcommand{\EG}{\mathcal{EG}}
\newcommand{\ELG}{\mathcal{ELG}}
\newcommand{\FG}{\mathcal{FG}}
\newcommand{\QP}{\mathcal{QP}}
\newcommand{\calCD}{\mathcal{C}}
\newcommand{\CO}{\mathcal{O}}

\newcommand{\la}{\lambda}

\newcommand{\bR}{\mathbb{R}}
\newcommand{\bC}{\mathbb{C}}

\newcommand{\bK}{\mathbb{K}}
\newcommand{\al}{\alpha}
\newcommand{\be}{\beta}
\newcommand{\ga}{\gamma}
\newcommand{\Si}{\Sigma}

\newcommand{\kP}{\mathbb{KP}}
\newcommand{\bP}{\mathbb{RP}}
\newcommand{\bcP}{\mathbb{CP}}
\newcommand{\brP}{\mathbb{RP}}
\newcommand{\gtn}{G_{2,n+1}}

\begin{document}

\title[Classifying real polynomial pencils]
{Classifying real polynomial pencils}

\author[J.~Borcea]{Julius Borcea}
\address{Department of Mathematics, Stockholm University, SE-106 91 Stockholm,
    Sweden}
\email{julius@matematik.su.se}
\author[B.~Shapiro]{Boris Shapiro}
\address{Department of Mathematics, Stockholm University, SE-106 91 Stockholm,
    Sweden}
\email{shapiro@matematik.su.se}

\subjclass{Primary 58K05; Secondary 12D10, 14P05, 26C10, 30C15}
\keywords{Real polynomial pencils, Grassmann discriminant, boundary-weighted
gardens.}

\begin{abstract}
Let $\bP^n$ be the space of all
homogeneous polynomials of degree  $n$ in two variables with real
coefficients. The standard discriminant
$\D_{n+1}\subset \bP^n$ is Whitney stratified according to the
number and the multiplicities of multiple real zeros. A real polynomial
pencil, that is, a line $L\subset \bP^n$ is called generic
if it intersects $\D_{n+1}$ transversally. Nongeneric pencils form
the Grassmann discriminant $\D_{2,n+1}\subset \gtn$, where
$\gtn$ is the Grassmannian of lines in $\bP^n$. We
enumerate the connected components of the set $\widetilde
\gtn=\gtn\setminus \D_{2,n+1}$ of all
generic lines in $\bP^n$ and relate this topic to the Hawaii conjecture and
the classical theorems of Obreschkoff and Hermite-Biehler.
\end{abstract}

\maketitle

\section{Introduction and main results}

In what follows by a
{\it pencil} $L=\{\al P + \be Q\}$ we will always mean a {\it real
polynomial pencil} of degree $n$ homogeneous polynomials in two real
variables, i.e., a real line in $\bP^n$ identified with the space of
all homogeneous degree $n$ real polynomials considered
up to a constant factor. Here $(\al : \be)\in \brP^1$ is a projective
parameter. In order to use derivatives
it will often be convenient to view homogeneous degree $n$ polynomials
in two variables as inhomeogeneous polynomials of degree at most $n$ in one
variable.  Any choice of a basis
    $(P,Q)$ in $L$ allows us to consider the real rational function $P/Q$;
    a different choice of basis leads to a rational function of the
    form $(AP+BQ)/(CP+DQ)$, which can
    be viewed as the postcomposition of the rational function $P/Q$ with
    the real linear fractional transformation $(Az+B)/(Cz+D)$ in the target
    $\bcP^1$. Thus all properties of real rational functions which are invariant
    under real linear fractional transformations in the target space are
naturally inherited by real polynomial pencils. For instance, the graph of a
real rational function $P/Q$ restricted to
    $\brP^1$ defines a finite branched covering $\brP^1\to \brP^1$. We
    call two rational functions $P_{1}/Q_{1}$ and $P_{2}/Q_{2}$ {\it
    graph-equivalent} if there exist diffeomorphisms of the source
    $\brP^1$ and the target $\brP^1$ sending the graph of $P_{1}/Q_{1}$
    to that of $P_{2}/Q_{2}$.  As a property which is invariant
under the
    postcomposition with a linear fractional transformation of the target space,
    the above graph-equivalence can be defined for the pencils
    $\al P_{1}+\be Q_{1}$ and $\al P_{2}+\be Q_{2}$.

\begin{figure}[!htb]
\centerline{\hbox{\epsfysize=2.8cm\epsfbox{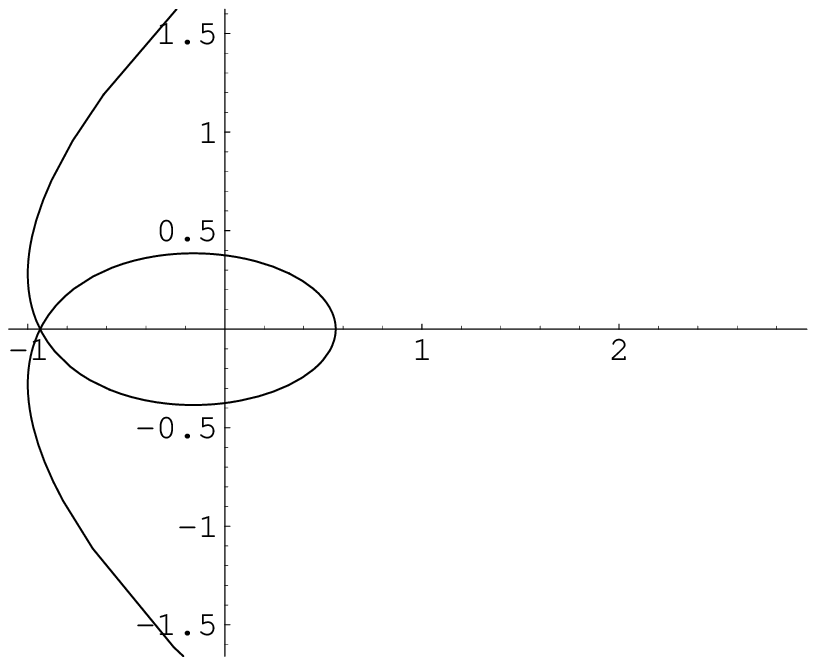}}
\hspace{0.5cm}\hbox{\epsfysize=2.8cm\epsfbox{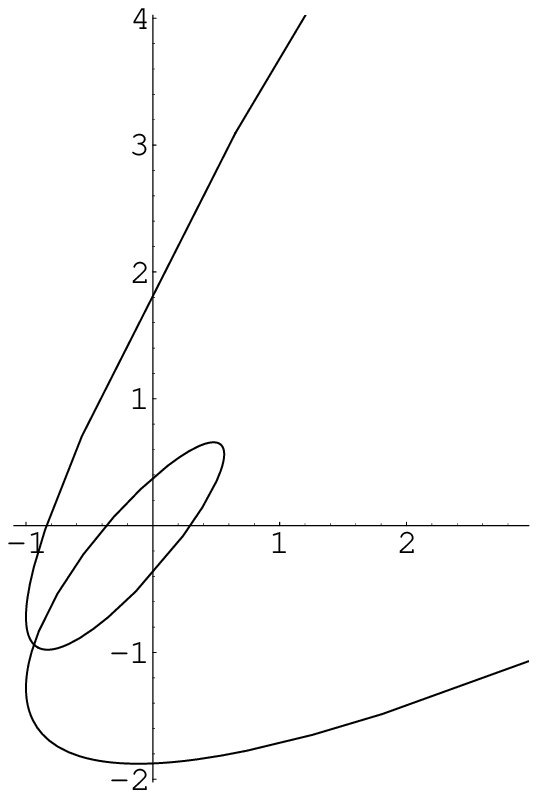}}
\hspace{0.5cm}\hbox{\epsfysize=2.8cm\epsfbox{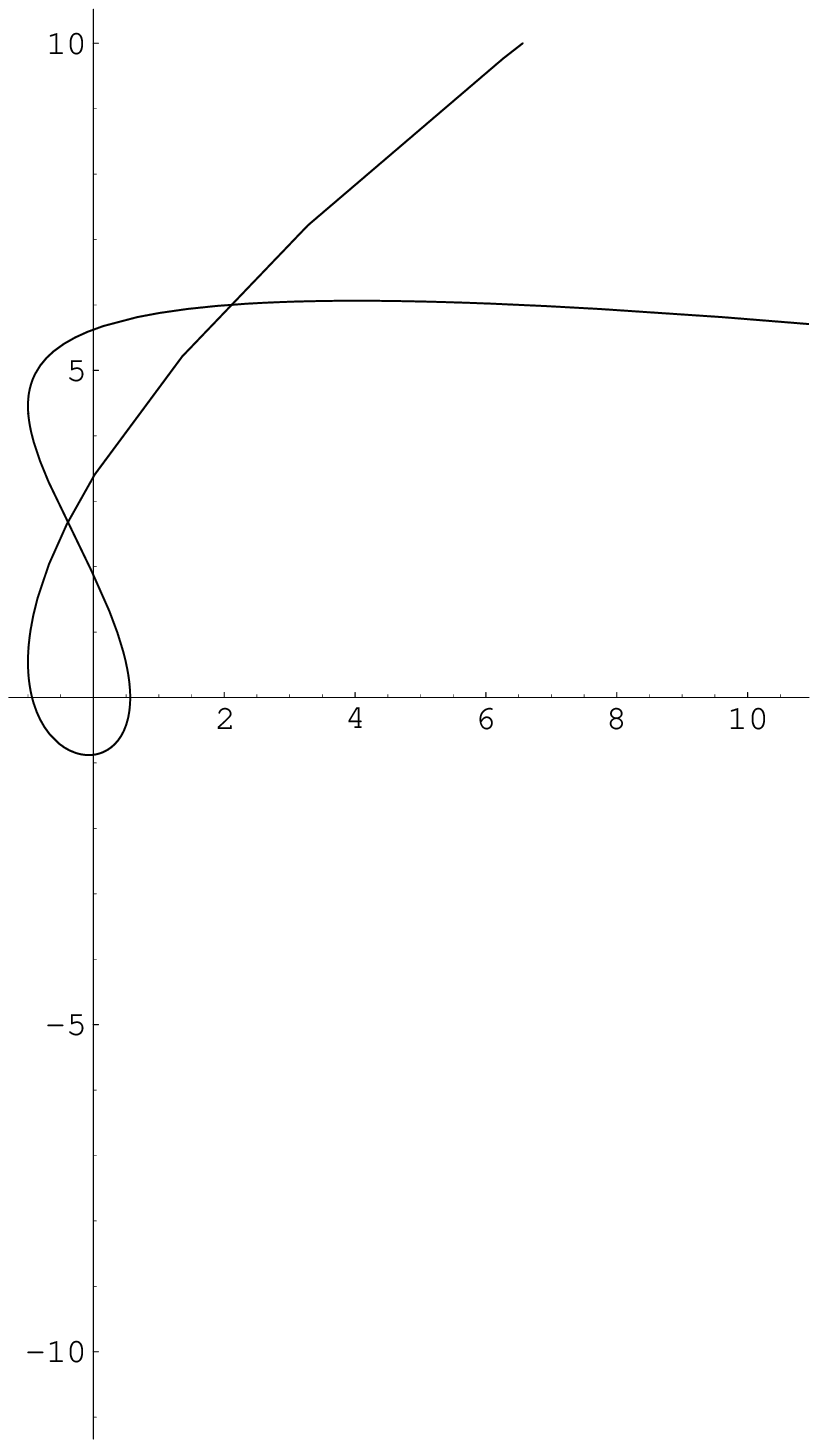}}}
\caption{Graph-equivalent and nonequivalent pencils}
\label{fig1}
\end{figure}

    The most classical notion of genericity for meromorphic
    functions/pencils requires that the function/pencil under consideration
    should have
    the maximal possible number of (simple) critical points with all
    distinct critical values. We will refer to this notion as {\it
    Hurwitz-genericity}, see \S 3. The classification of Hurwitz-generic
    real rational functions was carried out in details in \cite {NSV}.
    The violation of Hurwitz-genericity essentially occurs for two basic
reasons.
    Either several critical points collapse and form a degenerate
    critical point or some critical values collide but their corresponding
    critical
    points are still distinct. In the present paper we study a weaker
    notion of genericity than Hurwitz-genericity requiring only that all
    real critical points of the considered real rational
    functions/pencils stay simple, see Definition~\ref{d1} below. This
    notion is the natural
    counterpart of the absence of the collapse of critical points in the
    realm of real algebraic geometry. It still keeps some important
    information about the behavior of real rational functions/pencils and
    is closely related to the natural analog of the classical
    discriminant for the Grassmannian of two dimensional subspaces. One
    more important observation is that the violation of such genericity
    is detected in the source space instead of the target which is
    always more difficult. In short, we forbid singularities and do not
    care about multisingularities.






Our notion of genericity allows us in particular to give a complete solution
to the following problem.

\begin{problem}\label{pb1}
For which pencils  $L=\{\al P + \be Q\}$ is the number of
real zeros in this pencil constant, i.e., when is the number of real
solutions (counted with multiplicities) of the equation $\al P+\be Q=0$
independent of $\al/\be$?
\end{problem}

       An example of such a situation is provided by a well-known result
       of Obreschkoff, see \cite {Ob}, saying that  a  pencil
       $L=\{\al P+\be Q\}$ consists
of polynomials with only real (distinct) zeros if and only if both
$P$ and $Q$
have real (distinct) and interlacing zeros.
However, there exist pencils with a constant number of real zeros which are
not covered by Obreschkoff's result. For instance, one may consider the
pencil $L=\{\al P+\be P'\}$, where $P(x)=x^4+x^2-5x-4$.

\begin{figure}[!htb]
\centerline{\hbox{\epsfysize=2.8cm\epsfbox{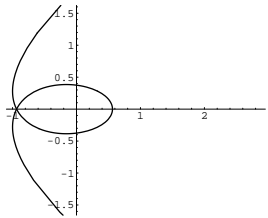}}}
\caption{When is the number of real zeros constant?}
\label{fig2}
\end{figure}

An easy observation is that a pencil $\{\al P + \be Q\}$  has a constant
number of real zeros if and only if the Wronskian $W(P,Q)=PQ'-QP'$
has
no real zeros, or in other words, $P$ and $Q$ form a fundamental
system for some second order linear ordinary differential equation.
Indeed, note that for any fixed $(\al:\be)$ the number of real zeros
of the equation
$\al P+\be Q=0$ equals the number of the real intersection points of
the rational
curve $\ga=(P,Q)$ with the line through the origin with slope
${\al}/{\be}$. If the number of real intersection points is constant then
there
should be no real tangent lines to $\ga$ passing through the origin.
But the points on $\ga$  where
the  tangent line to $\ga$  passes through the origin correspond
exactly to the real zeros of the Wronskian $W(P,Q)$,
see Figure~\ref{fig2}.

Let $\gtn$ denote as usual the Grassmannian of lines in $\bP^n$. The fact
that the behavior of the number of real zeros in a pencil
$L=\{\al P + \be Q\}$ is closely related to the properties of real
zeros of the Wronskian $W(P,Q)=PQ'-QP'$ justifies the following definition:

\begin{definition}\label{d1}
A real polynomial pencil $L=\{\al P + \be Q\}$ is called {\it
generic} if the Wronskian $W(P,Q)=PQ'-QP'$ has no multiple real
zeros and it is called {\it nongeneric} otherwise. The set $\D_{2,n+1}$ of all
nongeneric
real pencils in $\gtn$ is called the {\it Grassmann discriminant}.
\end{definition}

Clearly, the degree of the Wronskian of almost any pencil
in $\bP^n$ equals $2n-2$. If the Wronskian of
a pencil in $\bP^n$ is of the degree $2n-4$ or less then we consider this
pencil as degenerate (since its Wronskian has a double zero at $\infty$).

\begin{definition}\label{d2}
Two generic pencils are called {\it equivalent} if
they can be connected by a path through generic pencils, i.e., if they belong
to the same
connected component of the set $\widetilde
\gtn=\gtn\setminus \D_{2,n+1}$ of all generic pencils in $\bP^n$.
\end{definition}

Note that as defined above, the equivalence of two generic real pencils does
not necessarily imply their graph-equivalence since real critical values can
collide. The main question that we address below is the following.

\begin{problem}\label{pb2}
Enumerate  the equivalence classes of all generic pencils in $\bP^n$.
\end{problem}

The study of this topic originated from our attempt to solve the following
intriguing  conjecture of Craven, Csordas and Smith (cf.~\cite {CCS}; see also
\cite{ShS}).

\begin{conjecture}[Hawaiian Conjecture]\label{con1}
If a real polynomial $P$ has $2s$
nonreal zeros then the Wronskian
$W(P,P')=PP''-(P')^2$ has at most $2s$ real zeros.
\end{conjecture}

The main result of this paper -- Theorem~\ref{t1} below -- completely solves
Problem~\ref{pb2}. The answer
is given in terms of boundary-weighted gardens of total
weight $n$, a notion that we define and study in detail in \S 2 and \S 3.
The notion of garden of a real rational function provides also a natural
topological context for studying
Conjecture~\ref{con1} and related questions, see Conjectures~\ref{con2}
and~\ref{con3} in \S 5.

\begin{theorem}\label{t1}
The connected components in the space
     $\widetilde \gtn=\gtn\setminus \D_{2,n+1}$ of all
generic pencils in $\bP^n$ are in $1-1$ correspondence with the set
of equivalence classes of all boundary-weighted gardens of total
weight $n$.
\end{theorem}

 From Theorem~\ref{t1} and the arguments involving the Wronskian that we
mentioned earlier we immediately deduce the following answer to
Problem~\ref{pb1}.

\begin{figure}[!htb]
\centerline{\hbox{\epsfysize=10cm\epsfbox{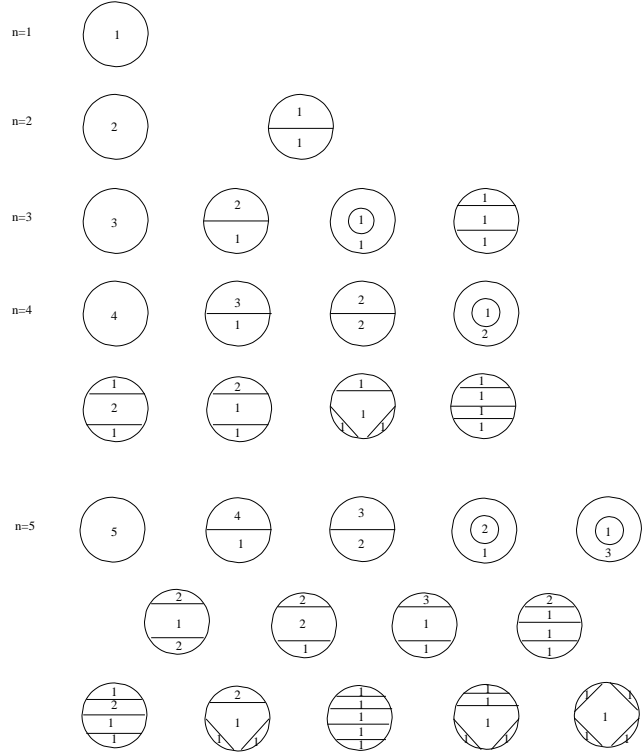}}}
\caption{Enumeration of equivalence classes for $n\le 5$}
\label{fig3}
\end{figure}

\begin{corollary}\label{cor1}
There exist $\left [\frac {n+1}{2}\right ]$
different components in
     $\widetilde \gtn=\gtn\setminus \D_{2,n+1}$ where the Wronskian
$W(P,Q)$ has no
     real zeros at all.
\end{corollary}

The values of the number of connected components for small values of
$n$ are $1, 2, 4, 8, 14, 28$ for $n$ equal to $1, 2, 3, 4, 5, 6$, respectively
(see Figure~\ref{fig3}). This sequence of integers was not recognized by the
online encyclopedia of integer sequences.

The structure of the paper is as follows. In \S 2 we define the notions of
garden, boundary-weighted garden and Morse perestroika and list some of
their properties. We further study these notions in \S 3, where we prove
the main results of the paper. In \S 4 we build on some of the aforementioned
ideas and obtain a simple new proof of a generalization of the famous
Hermite-Biehler theorem. Finally, \S 5 contains a number of conjectures and
open problems.

\medskip

\noindent
{\bf Acknowledgements.}  The authors are grateful to I.~Krasikov for
numerous enlightening discussions.
The second author is
obliged to A.~Eremenko, A.~Gabrielov, S.~Natanzon and A.~Vainshtein for
shaping his understanding of the topology of the space of real
rational functions and the Wronski map. The financial support and the
stimulating atmosphere of the program ``Topological aspects of real algebraic
geometry'' held in Spring 2004 at MSRI Berkeley are also highly appreciated.

\section{Preliminaries on gardens and gardening}

Let us first recall the notion
of garden of a real polynomial pencil as defined in \cite {EG1} and
\cite {NSV}.

\begin{definition}\label{d3}
The  {\it garden} $\CG(L)$ of a real
polynomial pencil $L=\{\al P(z) +\be  Q(z)\}$ is the set of all $z\in
\bcP^1$
for which the rational function $f_{L}=\frac {P(z)}{Q(z)}$ attains
real values.
\end{definition}

Note that the defining property of $\CG(L)$ is actually independent of the
choice of real basis $(P,Q)$ of the real polynomial pencil $L$.

\begin{figure}[!htb]
\centerline{\hbox{\epsfysize=3cm\epsfbox{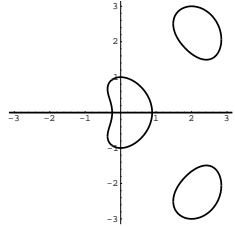}}}
\caption{Garden of the pencil $L=\{\al P+\be P'\}$, where
$P(x)=x^5-3x^4+10 x^3+10 x^2+9 x+13$}
\label{fig4}
\end{figure}

Let $z$ be an affine coordinate on $\bcP^1$. Observe that
$\CG(L)\subset \bcP^1$ is an algebraic
curve in the coordinates $(\Re{z}, \Im{z})$ which necessarily
contains $\bP^1\subset \bcP^1$ and is invariant under the complex
conjugation map $\tau: \bcP^1\to \bcP^1$. Note that the  singularities of any
garden occur exactly at the critical  points of $f_{L}=\frac
{P(z)}{Q(z)}$ where $f_{L}$ attains a real value. If such a critical
point has multiplicity $m\ge 2$ then at that point the garden has a
transversal intersection of $m$ nonsingular branches with angle
$\frac {\pi}{m}$ between any two neighboring branches. A critical
point with real critical value is called {\it simple} if its
multiplicity equals $2$.
A pencil $L=\{\al P(z) +\be  Q(z)\}$ is called {\it nonsingular}  if
the only critical points of $f_{L}$ with real critical values
    are real  and  simple. This implies that the only singularities of
its garden $\CG(L)$  are
transversal intersections of $\bP^1$ with other branches of $\CG(L)$.
The garden $\CG(L)$ of a nonsingular pencil $L$ will be called
{\it nonsingular} as well.
The aforementioned transversal intersections are called the {\it vertices} of
the garden.
Note that the vertices of $\CG(L)$ correspond exactly to the real
zeros of
$W(P,Q)$.

We need to describe nonsingular gardens in more details.
     A nonsingular garden $\CG(L)$  is the disjoint union of two  basic
parts
     $\CG(L)=\calCD(L)\cup \CO(l)$,  the {\it chord part} $\calCD(L)$ and the
(possibly empty)
{\it oval part} $\CO(L)$.
The chord part $\calCD(L)$ is the connected component of $\CG(L)$
containing $\brP^1$ while the oval part $\CO(L)$ is the complement
$\CG(L)\setminus\calCD(L)$. We call the edges connecting the
vertices of $\calCD(L)$ and not belonging to $\brP^1$ the {\it chords}.
The oval part $\CO(L)$
      consists of a number of $\tau$-invariant smooth closed curves
      called {\it ovals}. The connected components of $\bcP^1\setminus
\CG(L)$ are
      called the {\it faces} of the garden $\CG(L)$.
      Let us fix the
      standard metric on the image $\bcP^1$ such that the length of $\bP^1$
      equals $1$. If we  choose some basis $(P,Q)$ of the nonsingular pencil $L$
      under consideration then by using the rational
      function $f_{L}= \frac {P(z)}{Q(z)}: \bcP^1\to\bcP^1$ we can assign an
      extra piece of information
      to all elements of the garden $\CG(L)$.

\begin{definition}\label{ewg}
The {\it edge-weighted garden} $\EG(P/Q)$ of the rational function
$P/Q$ is the garden $\CG(L)$ of the pencil
      $L=\{\al P+ \be Q\}$ together with all edges, chords and
      ovals, each of these objects being endowed with the {\it weight} given by
the length of its respective image in the target
      $\bP^1$ under the rational function $P/Q$. By the {\it total weight} of
an edge-weighted garden  we mean the sum of the weights of all its edges,
      chords and ovals.
\end{definition}

\begin{remark}\label{r1}
Note that the image of an edge, chord or oval can
cover some interval of $\brP^1$ several times. The lengths/weights considered
in Definition~\ref{ewg} are total lengths obtained by counting multiplicities.
In particular, this implies that the total weight of the edge-weighted
      garden $\EG(P/Q)$
      equals  the degree of $P/Q$ as a map from $\bcP^1$ to $\bcP^1$.
\end{remark}

\begin{definition}\label{bwg}
A {\it boundary-weighted garden} is
a  nonsingular garden with positive integer weights assigned to each
boundary component of each face and satisfying the additional requirement
that $\tau$-symmetric faces are assigned equal weights. The
{\it total weight} of a boundary-weighted garden is the sum of the weights of
all boundary components contained in the closed upper hemisphere
other than ovals plus twice the weight of all ovals in the
closed upper hemisphere.
\end{definition}

\begin{figure}[!htb]
\centerline{\hbox{\epsfysize=4cm\epsfbox{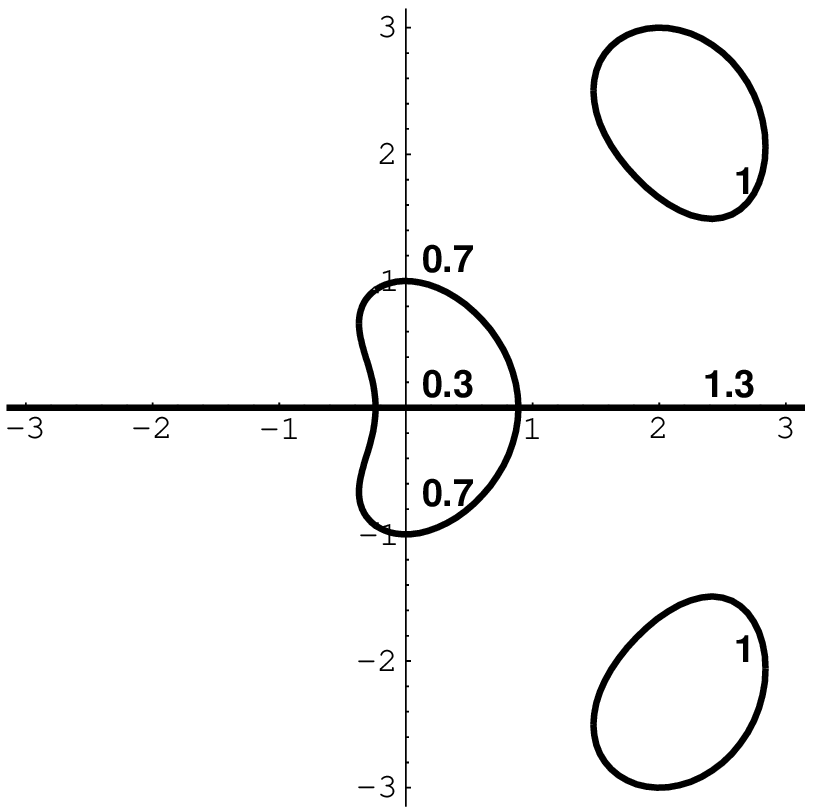}}
\hspace{1cm}\hbox{\epsfysize=4cm\epsfbox{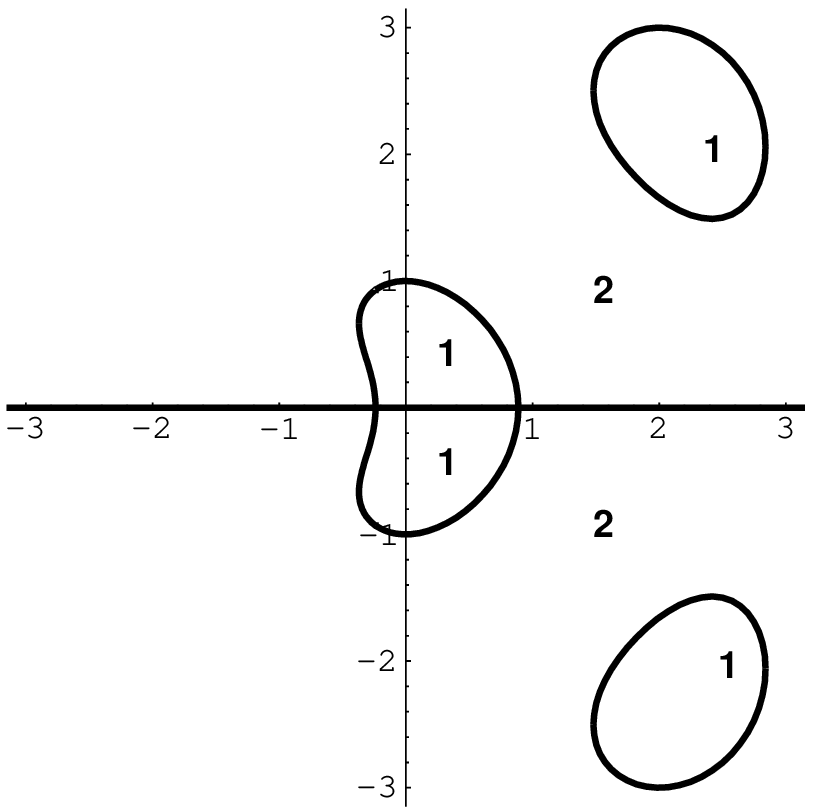}}}
\caption{An edge-weighted garden and its boundary-weighted image}
\label{fig5}
\end{figure}

      There is an obvious map $\la$ from edge-weighted
      gardens to boundary-weighted gardens obtained by assigning to each
    boundary component the sum of the weights of the elements contained
    in this boundary component. Note that the latter sum is either the sum of
the weights of all edges and chords if the
    boundary component contains them,
    or just the weight of an oval if the boundary component is an oval, see
Figure~\ref{fig5}. One can easily see that the image under $\la$ of an
edge-weighted garden $\EG(P/Q)$ is invariant under postcompositions of
the rational function $P/Q$ with real linear fractional transformations.
We may therefore associate to each nonsingular pencil a canonical
boundary-weighted garden in the following way.

\begin{definition}\label{cbwg}
The {\it boundary-weighted garden} $\FG(L)$ of a given nonsingular
pencil $L$ is the image under the map $\la$ of the edge-weighted garden
$\EG(P/Q)$, where $(P,Q)$ is some basis of $L$.
\end{definition}

Note that the integer placed in each face on Figure~\ref{fig5} is the weight of
the outer boundary component of the face (if the face is
multiconnected).
In order to describe connected components in the space of generic
real pencils we need to introduce the following equivalence relation
on the set of all boundary-weighted gardens of given total weight. In what
follows we will work with the half of a garden contained in the upper
hemisphere and assume that all operations are performed symmetrically.

\begin{definition}\label{d4}
By a {\it Morse perestroika} of a
boundary-weighted garden we mean
the following operation. Choose any face whose boundary contains
either two chords, two ovals
or a chord and an oval. Drag them together  and cut and paste them.
Under this operation two disjoint boundary
components will be glued together into one whose weight is the sum of
the weights of the former components. If the original face was simply
connected then it will be cut into two new faces. Its boundary will be
cut into two new boundary components whose weights are two arbitrary
positive integers which add up to the weigth of the former boundary.
Two boundary-weighted gardens that can be obtained from each
other by a sequence of Morse perestroikas are called {\it equivalent}.
\end{definition}

It is not difficult to see that the equivalence relation introduced in
Definition~\ref{d4} preserves the
total weight  and the number of chords of a garden.

\begin{figure}[!htb]
\centerline{\hbox{\epsfysize=2.8cm\epsfbox{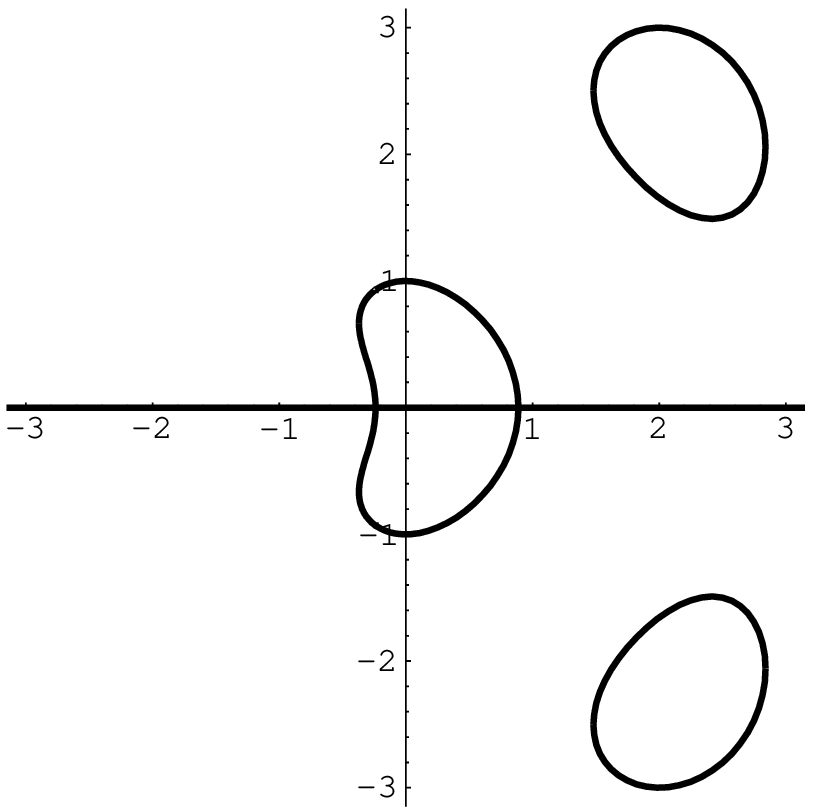}}
\hspace{0.5cm}\hbox{\epsfysize=2.8cm\epsfbox{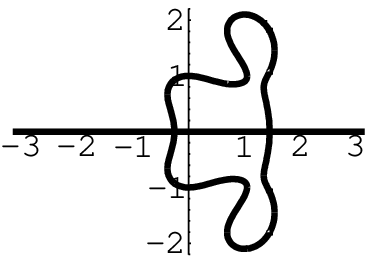}}
\hspace{0.5cm}\hbox{\epsfysize=2.8cm\epsfbox{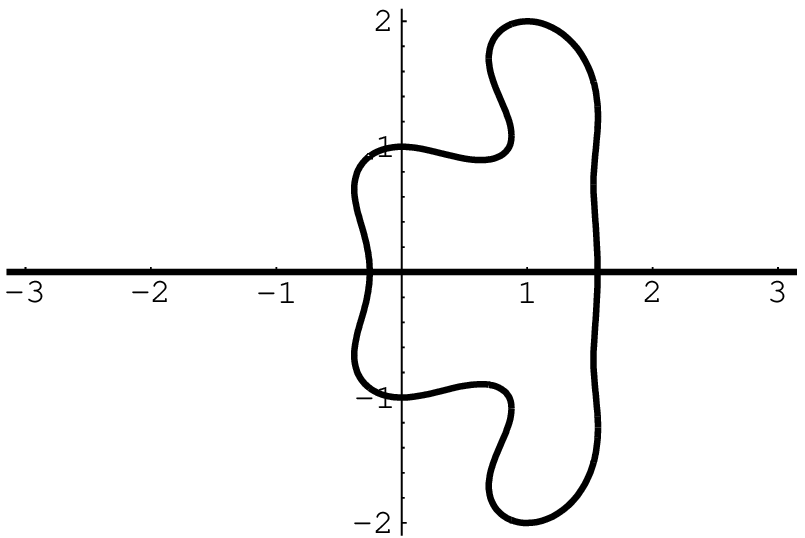}}}
\caption{Elementary perestroika in action}
\label{fig6}
\end{figure}

We have defined all the notions mentioned in Theorem~\ref{t1} and are now
ready to prove this theorem.

\section{Proofs}

We start with some generalities about $\D_{2,n+1}$ which can be
easily extended to linear polynomial families of higher dimensions.
The following important
mapping is called the {\it Wronski map}, see e.g.~\cite {EG2}. Let
$\bK$ denote $\bR$ or $\bC$.
     Introducing an affine coordinate $z$ on $\kP^1$ we can identify
     $\kP^n$ with the space of inhomogeneous polynomials of degree at
     most $n$ in the variable $z$.
     Consider now the map
$$W :\gtn\to \kP^{2n-2}$$ that sends a
2-dimensional linear polynomial subspace of $\kP^n$ to the linear span of
its Wronskian, i.e., the determinant of the $2\times 2$-matrix
$\left(\!\!\begin{array}{rr}
P(z) & Q(z)\\
P'(z) & Q'(z)
\end{array}\!\!\right)$,
where $(P(z),Q(z))$ is some basis of the chosen subspace.
Note that a change of basis in the given subspace amounts to multiplying the
Wronskian by a nonzero constant and that all such Wronskians are polynomials
in $z$ of degree at most $2n-2$.

      Several important facts are known about the map $W$.
      Over $\bC$ the map $W$ is finite and its degree  equals the degree of
$\gtn$
      under its Pl\"ucker embedding. The latter number equals the $n$-th Catalan
      number $\calCD_{n}=\frac 1 n
      \binom{2n-2}{n-1}$, see \cite {Go}. Moreover, the Wronski map is
perfectly adjusted
      to the Schubert cell decomposition of $\gtn$ constructed by using the
      natural complete flag in $\kP^n$ whose $i$-dimensional subspaces
      consist of all polynomials of degree at most $i$, where
$i=0,1,\ldots,n$. It turns out that over $\bC$ the degree of the
restriction of $W$ to
      any of the above Schubert cells equals the degree of this cell
      under the Pl\"ucker embedding of  $\gtn$, see \cite {EG2}.

Denote by $\D_{2n-2} \subset \kP^{2n-2}$ the standard
    discriminant in $\kP^{2n-2}$, that is, the set of all polynomials
    having a multiple zero over $\bK$. The Grassmann discriminant
    $\D_{2,n+1}$ introduced in Definition~\ref{d1} may alternatively be
characterized as follows.

\begin{definition}\label{d5}
The Grassmann discriminant
$\D_{2,n+1}\subset \gtn$ is  the inverse image $W^{-1}(\D_{2n-2})$
of $\D_{2n-2}$ under the Wronski map $W$.
\end{definition}

\begin{figure}[!htb]
\centerline{\hbox{\epsfysize=4cm\epsfbox{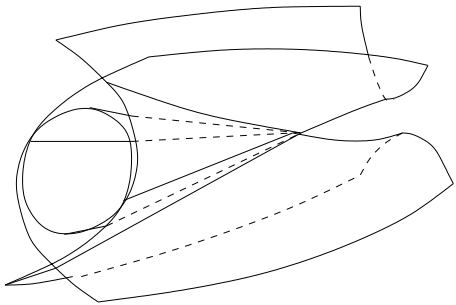}}}
\caption{Section of the $\D_{2,4}$-discriminant
transversal to the lift of $\rho$ to $G_{2,4}$ using tangent lines}
\label{fig7}
\end{figure}

\begin{lemma}\label{l1}
The discriminant $\D_{2,n+1}$ consists of two
irreducible components $U$ and $V$.  The first component $U$ is the
closure of the
set of all lines in $\kP^n$ tangent to $\D_{n}\subset \kP^n$
at its  smooth points. The second component $V$ is the set of all
lines passing through the stratum
$\Si_{3}\subset \D_{n}$, where $\Si_{3}$ consists of all polynomials
having a root over $\bK$ of multiplicity exceeding $2$ (compare with
\cite {GS} and see Figure~\ref{fig7}).
\end{lemma}

\begin{proof} Take a pencil $L=\{\al P +\be Q\}$ and consider the  matrix
$$M_{L}=\left(\!\!\begin{array}{rrr}
P(z)& P'(z)& P''(z)\\
Q(z)& Q'(z)& Q''(z)
\end{array}\!\!\right).$$
If the Wronskian
$W(P,Q)=\left|\!\begin{array}{cc}
P(z)& P'(z)\\
Q(z)& Q'(z)
\end{array}\!\right|$
has a multiple zero at some $z_{0}$ then
$$\left|\!\begin{array}{cc}
P(z_{0})& P'(z_{0})\\
Q(z_{0})& Q'(z_{0})
\end{array}\!\right|=
\left|\!\begin{array}{cc}
P(z_{0})& P''(z_{0})\\
Q(z_{0})& Q''(z_{0})
\end{array}\!\right|=0. $$

	     The latter conditions can be satisfied in two different
ways. Either
	     there exists $z_{0}$ such that
	     $P(z_{0})=Q(z_{0})=0,$ i.e.,  the first column in
	     $M_{L}$ vanishes  at $z_{0}$, or  the first column
	     never vanishes but there exists $z_{0}$ such that the first and
	     the second rows are linearly dependent.
	     The first situation corresponds to the case when the rational
	     curve $(P(z),Q(z))$ passes through the origin and the
	     corresponding pencil in $\kP^n$ is tangent to $\D_{n}$. The second
	     situation means that there exists a linear combination of $P$
	     and $Q$ which vanishes up to a cubic term, i.e., the pencil
	     intersects $\Si_{3}$, which means geometrically
	     that the curve $(P(z),Q(z))$ has a tangent line  at some
	     inflection point passing through the origin.
\end{proof}

For the sake of completeness let us present without proof yet
another characterization
of $\D_{2,n+1}$.
    The {\it standard rational normal curve}  $\rho: \kP^1\to
\kP^n$ is the curve consisting of all  degree $n$ polynomials with an
$n$-tuple root.
Given a complete projective flag $f$
in $\kP^n$ we associate to $f$ the standard Schubert cell
decomposition $\mathfrak{S}_{f}$ of $\gtn$ whose cells consist of all
$2$-dimensional
projective subspaces with a given set of dimensions of intersections
with the subspaces of $f$. The cells are labeled by Young
diagrams with at most
two rows of length not exceeding $n-1$.
Given a rational curve $\ga: \kP^1\to\kP^n$ one defines its
{\it flag lift} $\ga_{F}: \kP^1\to F_{n+1}$ to be the curve consisting of
all osculating flags to $\ga$. As is well known, the same definition applies
in fact to any projective algebraic curve.

\begin{proposition}\label{p1}
The component $U$ (respectively, $V$) of $\D_{2,n+1}$ is the union of
the Schubert cells
$\bigcup_{f\in \rho_{F}}C_{1,1}(f)$
(respectively, $\bigcup_{f\in \rho_{F}}C_{2,0}(f)$), where $f$ runs over the
flag lift $\rho_{F}$ of the standard  rational curve $\rho$. Here
$C_{1,1}(f)$ is the cell of codimension two in $\gtn$ whose Young
diagram with respect to $f$ is $(1,1)$ while $C_{2,0}$ is the cell whose
Young diagram with respect to $f$ is $(2,0)$.
\end{proposition}

In order to complete the proof of Theorem~\ref{t1} we need several additional
definitions and constructions. Let us first recall the following classical
definition.

\begin{definition}\label{d6}
A pencil $L=\{\al P(z) +\be Q(z)\}$ of degree $n$
polynomials
is called {\it Hurwitz-generic} if the rational function
$f_{L}=\frac {P(z)}{Q(z)}$ has $2n-2$ distinct critical points with distinct
critical values and it is called {\it Hurwitz-nongeneric} otherwise.
\end{definition}

\begin{remark}\label{r2}
As we already noted in the introduction, a real rational function of the
form $(AP+BQ)/(CP+DQ)$ may
    be viewed as the postcomposition of the rational function $P/Q$ with
    the linear fractional transformation $(Az+B)/(Cz+D)$ in the target
    $\bcP^1$. This shows that the property introduced in
Definition~\ref{d6} is independent of the choice of basis of the pencil $L$.
\end{remark}

\begin{definition}\label{d7}
The {\it Hurwitz discriminant} is the subset $\calH_{2,n+1}\subset \gtn$
consisting of all
Hurwitz-nongeneric pencils in the Grassmannian of lines in $\bP^n$.
\end{definition}

Clearly, any real Hurwitz-generic pencil $L$ is generic in the sense of
Definition~\ref{d1}. Moreover, such a pencil is also nonsingular, i.e., it
has a nonsingular garden $\CG(L)$. Indeed, any complex critical point
together with its complex conjugate form a pair that cannot have a real
critical value. This proves the following lemma.

\begin{lemma}\label{l2}
The Grassmann discriminant $\D_{2,n+1}$ is always contained in
the Hurwitz discriminant $\calH_{2,n+1}$.
\end{lemma}

We say that a nonsingular garden
    is {\it directed\/} if its edges, chords and ovals are directed in
such a way that
the boundary of each face becomes a directed cycle. This means that
any given face will lie either to the right of any of its boundary components
or to the left of any such component whenever we follow the direction that
has been assigned to a boundary component. The faces that
lie to the left of all of their boundary components are called
{\it positive\/} while faces lying to
    the right of their boundary components are called {\it negative\/}.
All neighbors of positive
    faces are negative and vice versa. Obviously, in order
to direct a garden it suffices to direct any one of its edges.
Therefore, there exist exactly two possible ways of directing a garden and
these are opposite to each other in the sense that the one is obtained from
the other by reversing the direction of every edge.

By a {\it proper (cyclic) labeling\/} of a directed
boundary-weighted garden with $2k$ vertices (equivalently, with $k$ chords)
we understand the labeling of its vertices by symbols
$1,\ldots, 2k$  satisfying the following condition: for each boundary component
the number of decreases (downs) between consecutive labels when we
traverse the labels of the vertices in the order
prescribed by the component's direction does not exceed the weight of this
component.
An {\it involution\/} of a properly labeled directed
boundary-weighted garden is an operation that reverses both its orientation
and the (cyclic) order of the labels by sending label $i$ to label $2k-i$.

Given a Hurwitz-generic rational function $f=P/Q$ and fixing the
orientation of the target $\brP^1$ one gets the  properly  directed
and labeled
boundary-weighted garden $\widetilde {\CG}(f)$ of $f$ by taking its
garden with the induced orientation of all elements and their induced weights
together with the cyclic labeling of
its vertices coming from the real critical values of $f$.
Note that for any other choice of basis of the pencil $L=\{\al P +\be Q\}$ the
resulting garden either coincides with $\widetilde {\CG}(f)$ or may be
obtained from $\widetilde {\CG}(f)$ by an involution.

An equivalent version of the following result was stated and proved by means
of rational functions in \cite  {NSV}.

\begin{theorem}\label{t2}
Let $\calH_{2,n+1}$ denote the divisor of all
Hurwitz-nongeneric pencils. The connected components in the space $\widetilde
H_{2,n+1}=\gtn\setminus \calH_{2,n+1}$ of all Hurwitz-generic pencils
are in $1-1$ correspondence with the set of all properly directed and
cyclicly labeled gardens of
weight $n$ modulo the action of the involution.
\end{theorem}

Recall from Lemma~\ref{l2} that the Grassmann discriminant $\D_{2,n+1}$ and
the Hurwitz discriminant $\calH_{2,n+1}$ satisfy
$\calH_{2,n+1}\supset \D_{2,n+1}$. For our further purposes we need the
following description of $\calH_{2,n+1}$.

\begin{theorem}\label{t3}
The Hurwitz discriminant $\calH_{2,n+1}$ is the union of four
real discriminants $U$, $V$, $W$ and $Z$, where $U$ and $V$ are defined in
Lemma~\ref{l1} and $W$ and $Z$ are two real algebraic hypersurfaces with the
same complexification, namely the hypersurface of all coinciding critical
values. More precisely, $W$ is the set of all
real pencils  $L=\{\al P +\be Q\}$ for which the rational
function $f_{L}=\frac
{P(z)}{Q(z)}$ has
two real critical points with coinciding real critical value, while
$Z$ is the set of all real pencils  $L=\{\al P +\be
Q\}$ for which the rational function $f_{L}=\frac {P(z)}{Q(z)}$ has
two complex conjugate critical points with coinciding (and therefore
real) critical value.
\end{theorem}

Our plan is as follows. We will show that by crossing $W$ one can realize any
admissible relabeling of a given cyclicly labeled boundary-weighted
garden and that by crossing $Z$ we can realize
any of its admissible Morse perestroikas. These two facts will be easy
corollaries of the following statements.

\begin{theorem}\label{t4}
Any edge-weighted garden $\CG$ of total weight $n$ is realized by a real
rational function of degree $n$. Moreover, the set of all
real rational functions
with a given egde-weighted oriented garden is path-connected.
\end{theorem}

Here by an edge-weighted garden of total weight $n$ we
understand an abstract embedded  $\tau$-invariant ``graph'' containing
$\brP^1$ with vertices only of
even multiplicity and possibly containing a number of
$\tau$-invariant ovals considered up to a diffeomorphism of the plane. All
edges, chords and ovals of this ``graph'' are equipped with positive
weights. Moreover, ovals have integer weights. Finally, for any boundary
component the sum of all weights in this component is a positive integer and
the sum of the weights of all elements in this ``graph'' equals $n$.
It is important to note that in Theorem~\ref{t4} we
{\em do not assume} that $\CG$ is a nonsingular garden and that
we actually allow arbitrary complex critical points with real critical values.

\begin{proof}[Proof of Theorem~\ref{t4}]
The proof is based on ideas similar to those used in the proof of Theorem 1
in \cite {NSV} and so we will only sketch it here. (The only major
difference compared to \cite {NSV} is that we allow singular
gardens.) We want to
construct a topological branched covering $ \bcP^1\to \bcP^1$ which is
invariant under complex conjugation and whose garden is isomorphic to $\CG$.
This will
prove the realization theorem,
since by Riemann's uniqueness theorem
there exists a unique complex structure on $\bcP^1$ for which this
topological covering is holomorphic. The orientation of the garden
uniquely specifies which
of its faces should be mapped to the upper hemisphere and also which faces
should be mapped to the lower hemisphere.
(Neighboring faces are always mapped to opposite hemispheres.)
Each open face of the garden is
a topological surface of genus $0$. The normalization of its closure
is a closed topological surface with boundary.
Now for any face of $\CG$ consider the total weight of its boundary
components, that is, the number of times each boundary component should
traverse $\brP^1$.
The Riemann-Hurwitz formula determines the number of simple complex critical
points the face under consideration should contain. We also know in which
hemisphere the corresponding critical values should lie. Let us now
recall some definitions from \cite[\S 3.1]{NSV}.
Denote by $\Lambda^+$ the upper
hemisphere $\{z\in \bcP^1 \; \mid \;\text{Im}\; z \ge 0 \}$ and by
$P$ a  genus $g$ topological surface with a boundary consisting of $k$
connected components. Consider the set $ \calH_{g,m}^k$ of all generic
degree $m$ branched coverings of the form $\phi : P\to \Lambda^+$ and let
$a_{1},\ldots,a_k$ be all the distinct connected components of $\partial
   P$. Given a partition $(m_1,...,m_k)\vdash m$
denote by $\calH_{g,m}^k(m_{1},\ldots,m_k)\subset \calH_{g,m}^k$
   the subset of maps
$\phi : P \to  \Lambda^+$ such that
$\deg \phi\vert_{a_{i}}=m_{i}$ for $i=1,\ldots , k$.
Obviously,
$$
\calH_{g,m}^k=\bigcup_{(m_1,...,m_k)\vdash m}\calH_{g,m}^k(m_{1},
\ldots,m_k).
$$

Let $f$ be a face in the upper hemisphere of $\bcP^1\setminus\CG$ and
consider the
space $\calH_{f}$ of all branched coverings from the normalization of
the closure of $f$ to $\Lambda^{\pm}$, where $\Lambda^{\pm}$ is the
upper or lower hemisphere depending on where $f$ should be mapped
according to the chosen orientation. Lemma 2 in \cite {NSV} shows that for
any partition $(m_1,...,m_k) \vdash m$ the space
$\calH_{g,m}^k(m_{1},\ldots,m_k)$ is path-connected. In particular, this
implies that each space $\calH_{f}$ is path-connected. We need the following
result.

\begin{lemma}\label{lDesc}
Let $\CG$ be an egde-weighted oriented garden and fix an arbitrary set of
(critical) values for the vertices belonging to its chords. Denote by
$Rat_{\CG}$
the set of all real rational functions with
      egde-weighted oriented garden $\CG$ and having these prescribed critical
values. Then $Rat_{\CG}$ is homeomorphic to
      $\Pi_{f\in Ind_{\CG}} \calH_{f}\times (\brP^1)^q,$
      where $q$ is the number of different connected components of
      $\CG$ containing vertices -- i.e., critical points with real
      critical values -- and $Ind_{\CG}$  is the index set of all faces $f$ in
the upper hemisphere of $\bcP^1\setminus\CG$.
\end{lemma}

\begin{remark}
Note that for a nonsingular garden with a positive
     number of vertices one has $q=1$ since all its vertices belong to
      $\brP^1$. However, the singular gardens considered in Theorem~\ref{t4}
might contain singular ovals with vertices which are not connected to
      $\brP^1$.
\end{remark}

\begin{proof}[Proof of Lemma~\ref{lDesc}]
Let us show first that by assigning all real critical
	values and picking an arbitrary map $\phi_{f}$ from each space
$\calH_{f}$
	for $f\in Ind_{\CG}$ we can glue together all the $\phi_{f}$'s into
precisely one half of a unique real
	rational function from $Rat_{\CG}$. This follows simply from
the fact that the real critical values determine exactly which parts of the
	boundary components of $\phi_{f_{i}}$ and $\phi_{f_{j}}$ for any two
neighboring faces $f_{i}$ and $f_{j}$ should be
	identified (glued together). Indeed, by gluing together all the
$\phi_{f}$'s for all $f\in Ind_{\CG}$ according to this recipe we get a
unique map from $\Lambda^{+}$ to $\bcP^1$. We may then take the
	conjugate copy of the latter map and glue the two halves
together along $\brP^1$ into a sphere
     $\bcP^1$, thus obtaining a unique final map $\bcP^1\to \bcP^1$. One can
easily see that the final map is the topological branched covering
     that satisfies all the properties required. It just remains to
     notice that in order to assign all real critical values for an
     edge-weighted garden it is necessary and sufficient to assign
     arbitrarily just one real critical value for each connected
     component of $\CG$ containing vertices. The
     critical values of the remaining vertices in each such component will
then be automatically
     restored from the set of weights of the chords and ovals in the respective
     component.
\end{proof}

To finish the proof of Theorem \ref{t4} just notice that the Cartesian
product of path-connected topological spaces is path-connected.
\end{proof}

\begin{corollary}\label{c2}
Any admissible  Morse perestroika of a given
nonsingular boundary-weighted garden is realizable.
\end{corollary}

\begin{proof}
Any singular garden that occurs while performing an arbitrary generic
perestroika contains just two simple complex conjugate critical points
with a common real critical value. It follows from Theorem~\ref{t4} that such a
garden can be realized by a rational function. Any small generic
1-parameter deformation of this rational function will necessarily produce
the required
perestroika. Indeed, in any such deformation the imaginary part of the
interesting critical value will necessarily change signs while the rest
of the garden will topologically stay the same.
\end{proof}

\begin{theorem}\label{t5}
The set of all real rational functions with a given
boundary-weighted oriented garden is path-connected.
\end{theorem}

\begin{proof} We use an argument similar to that of \cite {EG1}. Let $\CG$ be
an oriented  boundary-weighted garden and denote by $\ELG$ the
set of all possible edge-weighted gardens whose boundary-weighted
gardens coincide with $\CG$, see Figure~\ref{fig5}.
Enumerating arbitrarily all chords and edges in $\CG$
    and denoting the weight of the $i$-th chord by $w_{c,i}$, the
    weight of the $j$-th edge by $w_{e,j}$ and the
    weight of the $m$-th boundary component by $w(B_{m})$
    we get the following system of linear inequalities (one for each edge and
    chord) and linear equations  (one for each boundary component other than
an oval) satisfied by edge weights for all gardens in $\ELG$
\begin{equation}\label{e1}
\begin{cases}
       w_{c,i}>0, \\
       w_{e,j}>0, \\
       \sum_{i_{l}\in B_{m}}w_{c,i_{l}}+\sum_{j_{l}\in
	  B_{m}}w_{e,j_{l}}=w(B_{m}).
\end{cases}
\end{equation}

Let $Sol_{\CG}$ denote the set of all solutions to system~\eqref{e1}.
Obviously,
$Sol_{\CG}$ is a nonempty convex polytope. For any solution of~\eqref{e1}
we get an edge-weighted oriented garden. By Theorem~\ref{t4} the set
of all real rational functions realizing such a garden is path-connected.
Therefore, the set of real rational functions with a given oriented
boundary-weighted garden is actually fibered over a contractible base with
isomorphic path-connected fibers. (Note that by Lemma~\ref{lDesc}
the topology of the fiber does not depend on the particular weights of the
chords.) Thus the total space of fibration is path-connected.
\end{proof}

\begin{corollary}\label{c3}
Any admissible relabeling of a given
boundary-weighted and cyclicly labeled garden is realizable.
\end{corollary}

\begin{proof} Take any admissible labeling of a given boundary-weighted
garden. Place its labels arbitrarily on $\brP^1$ in an
order-preserving way, i.e., assign real critical values to all real
critical points. Then one can restore the weights of all the chords and
edges of the garden. These weights will necessarily satisfy
system~\eqref{e1} given above.  Having
done so for two different labelings and using the fact that the set
of rational functions in Theorem 5 is path-connected we conclude that
we can find a path from the first rational function to the second through
rational functions with the same boundary-weighted garden.
\end{proof}

Corollary~\ref{c3} completes the proof of Theorem~\ref {t1}. Let us
now use this theorem to deduce Corollary~\ref{cor1}.

\begin{proof}[Proof of Corollary~\ref{cor1}]
For a real polynomial pencil $L=\{\al P+ \be Q\}$ the
Wronskian $W(P,Q)$ has no real zeros if and only if the garden
$\CG(L)$ (as well as its equivalence class) has no chords. Among all
equivalence classes of boundary-weighted gardens of total weight $n$
there are exactly $\left[\frac{n+1}{2}\right]$ classes with no
chords. This is because the boundary-weighted garden of every such class
consists of $\brP^1$ and at most one
additional oval whose respective weights are integers $k$ and $l$
that satisfy $1\le k\le
n$, $0\le l$, $k\equiv n \bmod 2$ and $k+2l=n$ (cf.~Definition~\ref{bwg} and
Remark~\ref{r1}). The cases when $n\le 5$ are illustrated in Figure~\ref{fig3}.
\end{proof}

\section{Real pencils and the Hermite-Biehler theorem}

The properties of a real pencil $\{\al P +\be Q\}$ or, equivalently, of
the plane rational curve $\ga=(P,Q)$ are also involved in the following
well-known result.
The classical Hermite-Biehler theorem asserts that given two
polynomials $P$ and $Q$ with real coefficients and of degrees $n$
and $n-1$, respectively, the zeros of the complex polynomial $S=P+iQ$ have
(nonzero) imaginary parts of the same sign if and only if
$P$ and $Q$ have real distinct and interlacing zeros.
In fact, if $\mu$ is an arbitrary complex number and $\sharp_{+}$ denotes
the number of zeros of the polynomial $S_{\mu}:=P+\mu Q$ lying in
the upper half-plane then the following more general result is known to be
true, see \cite{Ga}.

\begin{proposition}\label{p2}
In the above notation
consider  the plane real rational curve $\ga_{\mu}$ given by
$(P+\Re \mu \cdot Q,\Im \mu \cdot  Q)$. Then $\sharp_{+}$ equals the
winding number of
$\ga_{\mu}$ around the origin.
\end{proposition}

Below we give a new proof of the
generalized Hermite-Biehler theorem for all pairs $(P,Q)$ of real
polynomials. In particular, our method yields a simple proof of the
main result in \cite {HDB}.

\begin{proposition}\label{p3}
For given polynomials $P$ and $Q$ with real
coefficients the complex polynomial $S_{\mu}=P+\mu Q$ with $\mu\notin
\bR$ has a real zero if
and only if $P$ and $Q$ have a common real zero.
\end{proposition}

\begin{proof} Indeed, if $P(x)$ and $Q(x)$ have a common real zero
$x_{0}$ then $S_{\mu}(x_{0})=0$. On the other hand, if for some
$x_{0}\in \bR$ one has $S_{\mu}(x_{0})=0$ then $P(x_{0})+\Re \mu \cdot
Q(x_{0})=0$
and $\Im \mu \cdot Q(x_{0})=0$, which immediately imply
$P(x_{0})=Q(x_{0})=0$ since $\Im \mu\neq 0$.
\end{proof}

A convenient geometric reformulation of this statement is as follows.
Denote by $Pol_{n}$ the space of all monic
degree $n$ polynomials with complex coefficients of the form
$S(z)=z^n+a_{1}z^{n-1}+a_{2}z^{n-2}+\ldots+a_{n}$
and let $\RD\subset Pol_{n}$ be the
hypersurface of all polynomials $S$ that have at least one real zero.
Finally, let $Res\subseteq
Pol_{n}$ be the hypersurface of all $S=P+iQ$ such that $P$ and $Q$
have a real common zero. In the literature on singularities $Res$ is
often called the (generalized) Whitney umbrella.

\begin{corollary}\label{c4}
The discriminant $\RD$ coincides with the resultant hypersurface $Res$.
\end{corollary}

\begin{remark}\label{r4}
In the definition of $Res$
we disregard the subvariety of real codimension two  where $P$
and $Q$ have common complex zeros.
\end{remark}

Given an arrangement of black and white
distinct points on $\bR$ we define its {\it canonical reduction} to be the
interlacing (possibly empty) arrangement obtained in the following way: if
our arrangement contains a pair of neighboring points of the same
color
then we remove these points and we continue this procedure until no such
removals can be performed. Note that the resulting canonical reduction
depends only on the initial (relative) order of the points in the given
arrangement and not on their exact locations on $\bR$.

\begin{corollary}[cf.~\cite{KS}]\label{c5}
The number of connected components in
$Pol_{n}\setminus \RD$ equals $n+1$ and these components can be
labeled by the canonical reductions as follows. Let
$P$ and $Q$ be polynomials with real coefficients of degree $n$ and $n-1$,
respectively. Assume that they have no common real zeros and that the
leading coefficient of $Q$ is positive. Then for the polynomial
$S_{\mu}=P+\mu Q$ with $\Im \mu \neq 0$ one has
$\sharp_{+}-\sharp_{-}=\kappa  T$, where $\sharp_{+}$
(respectively, $\sharp_{-}$) is the number of zeros of
$S_{\mu}$ in the upper (respectively, lower) half-plane, $\kappa$ is the sign
of $\Im \mu $ and $T$ is the number of zeros of $P$ appearing in the
canonical reduction of the real zeros of $P$ and $Q$.
\end{corollary}

\section{Final remarks}

As we already mentioned
in the introduction, the notion of garden of a real rational function
provides a natural topological
framework for investigating the Hawaii conjecture (Conjecture~\ref{con1}).
Indeed, given a polynomial $P$ of degree $n$ with real coefficients let us
consider the garden $\CG_{P}$ of the rational function $P'/P$
(cf.~Definitions~\ref{d3}--\ref{ewg}). Obviously, all zeros of $P$ lie on
$\CG_{P}$. We make the following conjecture.

\begin{conjecture}\label{con2}
Each chord of $\CG_{P}$ contains at least one nonreal zero of $P$.
\end{conjecture}

Note that Conjecture~\ref{con2} would immediately imply the Hawaii
conjecture since the real critical points of $P'/P$ are the same as the real
zeros of the Wronskian $W(P,P')$ and the latter are precisely the endpoints
of the chords in $\CG_{P}$ (cf.~\S 2).

It is natural to ask whether the Hawaii conjecture extends to classes of
rational functions other than logarithmic derivatives. Let $n$ be a positive
integer and denote by $\QP_n$ the
set of all nonidentically vanishing rational functions of the form
\begin{equation}\label{e2}
f(x)=\sum_{i=1}^{n}c_{i}P^{\al}_i(x),
\end{equation}
where $c_i\in \bR$ for $1\le i\le n$,
$\al$ is a real number satisfying $\al\le -1$
and $P_1,\ldots,P_n$ are second degree monic polynomials with real
coefficients without real zeros. Based on extensive numerical experiments,
we propose the following analog of Conjecture~\ref{con1} for the class
$\QP_n$.

\begin{conjecture}\label{con3}
If $f\in \QP_n$ then $f$ has at most $2n-1$ real critical points. Moreover,
if $\al$ is a negative integer then the following analog of
Conjecture~\ref{con2} holds: each chord of the garden $\CG_f$ of the real
rational function $f$ contains at least one nonreal zero of the polynomial
$\prod_{i=1}^{n}P_i(x)$.
\end{conjecture}

A possible way to attack Conjecture~\ref{con3} might be as follows. Let us
first recall the definition of a Tchebycheff system as given in 
e.g.~\cite {Ka}.

\begin{definition}
A linear $n$-dimensional space $V$ of smooth
real-valued functions defined on some interval $(a,b)$ ($a$ might be
equal to $-\infty$ and $b$ to $+\infty$) is called  a
{\it Tchebycheff system} if any nonidentically vanishing function $f\in V$
has at most $n-1$ real zeros on $(a,b)$ counted with multiplicities.
\end{definition}

\begin{problem}\label{pb7}
Let $P_1,\ldots,P_n$ be as in~\eqref{e2} and $\al\in \bR$ with $\al \le -1$.
Is it possible to extend the $n$-tuple of functions $(P^\al_{1}(x))',\ldots,
(P^\al_{n}(x))'$ to a Tchebycheff system  of dimension $2n$ on 
$(-\infty,+\infty)$?
\end{problem}

Note that an affirmative answer to Problem~\ref{pb7} would automatically
confirm the validity of Conjecture~\ref{con3}.

The main question about the classification of generic
pencils (Problem~\ref{pb2}) extends straightforwardly to polynomial families
with more than one parameter. However, a solution to the problem of
   enumerating connected components in other Grassmannians similar to
Theorem~\ref{t1} would first require an appropriate definition of the notion
of garden in these cases.

To conclude, let us formulate some related questions.


\begin{problem}\label{pb3}
What can one say about the topology of the space of generic pencils
$\widetilde\gtn=\gtn\setminus \D_{2,n+1}$? For instance, are connected
components in $\widetilde \gtn$
contractible? Note that this is true for polynomials without
multiple real roots.
\end{problem}

A real rational function of degree $n$ is called an {\it $M$-function}
if all its $2n-2$ critical points and critical values are real and
distinct. Any $M$-function of degree $n$ induces a degree $n$ map
$\brP^1\to \brP^1$ with exactly $2n-2$ branching points.

\begin{problem}\label{pb4}
What type of maps $\brP^1\to \brP^1$  of degree $n$
with $2n-2$ branching points can occur from $M$-functions?
\end{problem}

More precisely, given a map $\brP^1\to \brP^1$  of degree $n$
with $2n-2$ simple branching points let us label its $n$ real critical points
and the corresponding $n$ critical values cyclicly. Then we can associate to
this map the unique cyclic permutation of length $2n-2$ sending each critical
point to its critical value. Problem~\ref{pb4} may therefore be reformulated
as follows.

\begin{problem}\label{pb5}
What cyclic permutation can an $M$-function have?
\end{problem}

Note that Problem~\ref{pb5} is actually asking for a description
of all
possible shapes of the graphs of real rational $M$-functions -- a topic which
is standardly considered in
elementary calculus courses if one omits ``$M$-'' in the above
formulation. However, in the general case the answer to Problem~\ref{pb5}
seems to be unknown and quite nontrivial.

\begin{problem}\label{pb6}
Enumerate the connected components in the space of
real rational functions having only simple real critical points with
distinct critical values.
\end{problem}

The arguments given in the introduction show that all real pencils are
necessarily graph-equivalent in each such component. This project is the
intermediate situation between the one covered in \cite {NSV} and the one
described in the present article.


\begin{thebibliography}{99999}


\bibitem[CCS]{CCS}
T.~Craven, G.~Csordas, W.~Smith, {\em The zeros of derivatives of entire
functions and P\'olya-Wiman conjecture}, Ann. of Math. (2) {\bf 125} (1987),
405--431.


\bibitem[EG1]{EG1}
A.~Eremenko, A.~Gabrielov, {\em Rational functions with real critical points
and the B. and M. Shapiro conjecture in real enumerative geometry},
Ann. of Math. (2) {\bf 155} (2002), 105--129.

\bibitem[EG2]{EG2}
A.~Eremenko, A.~Gabrielov, {\em Degrees of real Wronski maps},
Discrete Comput. Geom. {\bf 28} (2002), 331--347.

\bibitem[Fu]{Fu}
M.~Fujiwara, {\em Einige Bemerkungen \"uber die elementare Theorie der
algebraischen Gleichungen}, T\^ohoku Math. J. {\bf 6} (1916), 102--108.

\bibitem[Ga]{Ga}
F.~R.~Gantmacher, The theory of matrices, Chelsea, NY, 1959.

\bibitem[Go]{Go}
L.~Goldberg, {\em Catalan numbers and branched coverings by the Riemann
sphere}, Adv. Math. {\bf 85} (1991), 129--144.

\bibitem[GS]{GS}
A.~Gorodentsev, B.~Shapiro, {\em On associated discriminants for polynomials
of one variable}, Beitr\"age Algebra Geom. {\bf 39} (1998), 53--74.

\bibitem[HDB]{HDB}
M-T.~Ho, A.~Datta, S.~P.~Bhattacharyya, {\em Generalizations of the
Hermite-Biehler theorem: the complex case}, Linear Algebra Appl. {\bf 320}
(2000),  23--36.

\bibitem[Ka]{Ka}
S.~Karlin, Total positivity, Vol.~I, Stanford Univ.~Press, Stanford,
CA, 1968.

\bibitem[KS]{KS}
B.~Khesin, B.~Shapiro, {\em Swallowtails and Whitney umbrellas are
homeomorphic}, J. Algebraic Geom. {\bf 1} (1992), 549--560.




\bibitem[NSV]{NSV}
S.~Natanzon, B.~Shapiro, A.~Vainshtein, {\em Topological classification of
generic real rational functions}, J. Knot Theory Ramifications {\bf 11}
(2002), 1063--1075.

\bibitem[Ob]{Ob}
N.~Obreschkoff, Verteilung und Berechnung der Nullstellen reeller Polynome,
VEB Deutscher Verlag der Wissenschafter, 1963.

\bibitem[ShS]{ShS}
T.~Sheil-Small, Complex polynomials, Cambridge Studies in Adv. Math. Vol.
{\bf 75}, Cambridge Univ. Press, Cambridge, UK, 2002.


\end{thebibliography}
\end{document}